\documentclass[12pt]{amsart}

\usepackage{tensor}     % for tensor indices

\usepackage{tikz}
\usetikzlibrary{calc, decorations.pathreplacing, arrows.meta, patterns}

\usepackage[colorlinks=true,urlcolor=blue, citecolor=red,linkcolor=blue,linktocpage,pdfpagelabels, bookmarksnumbered,bookmarksopen]{hyperref}
\usepackage[hyperpageref]{backref}
\usepackage{cleveref}
\usepackage{xcolor}
\usepackage{amsthm} %for citing inside theorem header
\usepackage{latexsym,amsmath,amssymb}
\usepackage{accents}
\usepackage[colorinlistoftodos,prependcaption,textsize=tiny]{todonotes}
\usepackage{a4wide}
\usepackage{soul}
\usepackage{mathtools} %For \chi with a much lower index (\mychi)
\usepackage{xparse} %For a new definiton of \chi with a slightly lower index
\usepackage{enumitem}		% for using [nosep] with itemize and enumerate

\usepackage{calc}

\usepackage{accents}

\definecolor{indigo}{rgb}{0.29, 0.0, 0.51}
\definecolor{p1}{gray}{0.4}
\definecolor{p2}{gray}{0.6}
\definecolor{p3}{gray}{0.98}
\definecolor{p4}{gray}{0.8}
\definecolor{p5}{gray}{0.9}

plus 6pt minus 12pt
\def\vp{\varphi}

\newcommand{\m}{\mathcal{M}}

\def\B{{B}}

\def\n{{\mathcal N}}
\def\cM{{\mathcal M}}

\def\n{{\mathcal{M}}}
\def\n{{\mathcal N}}

\def\S{{\mathbb S}}

\renewcommand{\div}{{\rm div}}

\newtheorem{theorem}{Theorem}
\newtheorem{lemma}[theorem]{Lemma}
\newtheorem{corollary}[theorem]{Corollary}

\newtheorem{remark}[theorem]{Remark}
\newtheorem{definition}[theorem]{Definition}

\def\x{\tilde x}

\newcommand{\dif}{\,\mathrm{d}}

\newcommand{\dd}{\dif}

\newcommand{\R}{\mathbb{R}}

\newcommand{\cN}{\mathcal{N}}

\newcommand{\pl}{\partial}

\newcommand{\brac}[1]{\left (#1 \right )}
\newcommand{\abs}[1]{\left |#1 \right |}

\newcommand{\la}{\mathopen{}\mathclose\bgroup\left\langle}
\newcommand{\ra}{\aftergroup\egroup\right\rangle}

\DeclareMathOperator{\Ric}{Ric}

\newcommand{\barint}{
\rule[.036in]{.12in}{.009in}\kern-.16in \displaystyle\int }

\newcommand{\barcal}{\mbox{$ \rule[.036in]{.11in}{.007in}\kern-.128in\int $}}

\def\mvint_#1{\mathchoice
          {\mathop{\vrule width 6pt height 3 pt depth -2.5pt
                  \kern -8pt \intop}\nolimits_{\kern -3pt #1}}%
          {\mathop{\vrule width 5pt height 3 pt depth -2.6pt
                  \kern -6pt \intop}\nolimits_{#1}}%
          {\mathop{\vrule width 5pt height 3 pt depth -2.6pt
                  \kern -6pt \intop}\nolimits_{#1}}%
          {\mathop{\vrule width 5pt height 3 pt depth -2.6pt
                  \kern -6pt \intop}\nolimits_{#1}}}

\numberwithin{theorem}{section} \numberwithin{equation}{section}

\def\XXint#1#2#3{{\setbox0=\hbox{$#1{#2#3}{\int}$}
     \vcenter{\hbox{$#2#3$}}\kern-.5\wd0}}

\let\latexchi\chi
\makeatletter
\renewcommand\chi{\@ifnextchar_\sub@chi\latexchi}
\newcommand{\sub@chi}[2]{% #1 is _, #2 is the subscript
  \@ifnextchar^{\subsup@chi{#2}}{\latexchi^{}_{#2}}%
}
\newcommand{\subsup@chi}[3]{% #1 is the subscript, #2 is ^, #3 is the superscript
  \latexchi_{#1}^{#3}%
}
\makeatother
\DeclareMathOperator{\stab}{Stab}
\DeclareMathOperator{\boch}{Boch}

\usepackage{mathtools} % MM: it's already in macros, but for some reason doesn't work

\title[Kato inequality]{Vectorial Kato inequality for $p$-harmonic maps\\ with optimal constant}

\author{Andreas Gastel}
\address[Andreas Gastel]{Fakult\"{a}t f\"{u}r Mathematik, Universit\"{a}t
Duisburg-Essen, Thea-Leymann-Str.~9, 45127 Essen, Germany
}
\email{andreas.gastel@uni-due.de}

\author{Katarzyna Mazowiecka}
\address[Katarzyna Mazowiecka]{
Institute of Mathematics, %
University of Warsaw,
Banacha 2,
02-097 Warszawa, Poland}
\email{k.mazowiecka@mimuw.edu.pl}

\author{Micha\l{} Mi\'{s}kiewicz}
\address[Micha\l{} Mi\'{s}kiewicz]{
Institute of Mathematics, %
University of Warsaw,
Banacha 2,
02-097 Warszawa, Poland}
\email{m.miskiewicz@mimuw.edu.pl}

\begin{document}

 \begin{abstract}
 We derive the sharp vectorial Kato inequality for $p$-harmonic mappings. Surprisingly, the optimal constant differs from the one obtained for scalar valued $p$-harmonic functions by Chang, Chen, and Wei \cite{ChaCheWei16}. As an application we demonstrate how this inequality can be used in the study of regularity of $p$-harmonic maps. Furthermore, in the case of $p$-harmonic maps from $B^3$ to $\S^3$, we enhance the known range of $p$ values for which regularity is achieved. Specifically, we establish that for $p \in [2, 2.642]$, minimizing $p$-harmonic maps must be regular.
 \end{abstract}

% \keywords{Minimizing fractional harmonic maps, homotopy theory, regularity theory, existence}
\sloppy

\maketitle
\tableofcontents
\sloppy
\section{Introduction}
For differentiable vector fields, or sections of a vector bundle with metric
connection $\nabla$, the trivial
pointwise inequality $|\nabla V|\ge|\nabla|V||$ is known as
Kato's inequality. It has been known for a long time that the constant
$1$ in this inequality can be improved for solutions of (systems of)
partial differential equations. This dates back at least as early as 1975
when Schoen, Simon, and Yau \cite{SSY} proved $|\nabla A|^2\ge\frac{n+2}n\,|\nabla|A||^2$
for the second fundamental form~$A$ of any minimal hypersurface in $\R^{n+1}$. This may not have been stated explicitly in \cite{SSY}, see \cite[Lemme 3]{Berard} for this interpretation.

For harmonic maps $u\colon \mathcal M^n\to \mathcal N^d$, the vector field to be considered is
$\nabla u$, giving estimates between $|\nabla^2 u|$ and $|\nabla|\nabla u||$.
The most striking application of Kato inequalities here is Okayasu's paper \cite{Okayasu}. He proved
$|\nabla|\nabla u||^2 \le \frac{n-1}n \,|\nabla^2 u|^2$ in the sense of
distributions, with optimal constant, and used
this to improve the results by Schoen and Uhlenbeck \cite{SU3} about partial
regularity for harmonic maps to spheres, resulting in optimal regularity
for maps to $\S^4$ and $\S^5$. Note, however, that Schoen and Uhlenbeck
already had some improved (but non-optimal) Kato-type inequality with constant
$\frac{2nd}{2nd+1}$ instead of $\frac{n-1}n$.

Improved Kato type inequalities are most significant in regularity proofs where
the values of constants are crucial. This is the case in all of the papers
mentioned so far, and more generally in many regularity proofs that perform
``dimension reduction of the singular set'' by ruling out conical solutions
--- representing the blow-up of a solution near a singular point --- by playing
out a stability inequality against some Bochner--Weitzenb\"ock formula.
As this technique has widespread applications, the search for optimal
constants in Kato's inequality is relevant for solutions of many types
of equations and systems.

For harmonic maps to targets of any dimension, Okayasu's constant $\frac{n-1}n$
is optimal. The situation is different for $p$-harmonic maps, which are defined as critical points of
\begin{equation}\label{def:Ep-energy}
 E_p(u) = \int_{\mathcal M} |\nabla u|^p \qquad \text{among $u\in W^{1,p}(\mathcal M, \mathcal N)$.}
\end{equation}
It is remarkable that the optimal constant for $p$-harmonic maps  is not the same in the scalar and the vectorial case. This may be one of the reasons why Kato-type inequalities for $p$-harmonic
maps have not been proven but recently. Let us write the inequalities in
the form $|\nabla^2u|^2\ge\kappa|\nabla|\nabla u||^2$ here.
In 2016, Chang, Chen, and Wei \cite{ChaCheWei16} proved a Kato inequality for $p$-harmonic
{\em functions} ($\R$-valued) with the optimal constant
$\kappa=\min\{1+\frac{(p-1)^2}{n-1},2\}$. In 2019, the first-named author \cite{Gastel19} studied a micropolar elasticity model which couples $p$-harmonic
maps with some other equation. In the course of a proof that the theory
is sort of well-behaved also for exponents $p$ slightly larger than $2$, he provided a first Kato-type inequality for targets other than $\R$, with some non-optimal $\kappa>1$ for suitable
combinations of $n$ and $p$. This prompted the
second- and third-named authors \cite{MM} in 2023 to study the $n=2$ case more
thoroughly, with the optimal result $\kappa=\min\{1+(p-1)^2,2\}$, which happens to
coincide with the constant for scalar-valued case mentioned above.

These results, combined with the harmonic case $p=2$, might lead one to
guess that $\kappa=\min\{1+\frac{(p-1)^2}{n-1}, 2\}$ could be the
best constant in the general case $n\ge2$, $p>1$. But \cite{MM} also gave numerical evidence that in some cases the constant should be worse than that. This was our motivation to tackle the problem of the optimal Kato constant for
$p$-harmonic maps in general. Our result is that the constant is given by
\[
    \kappa(p,n) =\left\{ \begin{array}{lll}
                     1+\frac{(p-1)^2}{n-1} &\quad  \text{for }&1\le p\le1 + \frac{n-1}{\sqrt{2n}-1}\\
                        2-\frac{(p-\sqrt{2n})^2}{(\sqrt{n}-\sqrt{2})^2} &\quad \text{for }&1 + \frac{n-1}{\sqrt{2n}-1}\le p \le \sqrt{2n}\\
                     2 &\quad \text{for }& p\ge \sqrt{2n}.
                    \end{array}
\right.
\]
The constant above is sharp in a pointwise sense: it is achieved by some map that solves the $p$-harmonic equation at the given point. 
This is the main result of the present paper and it is proved in \Cref{ch:Kato}.

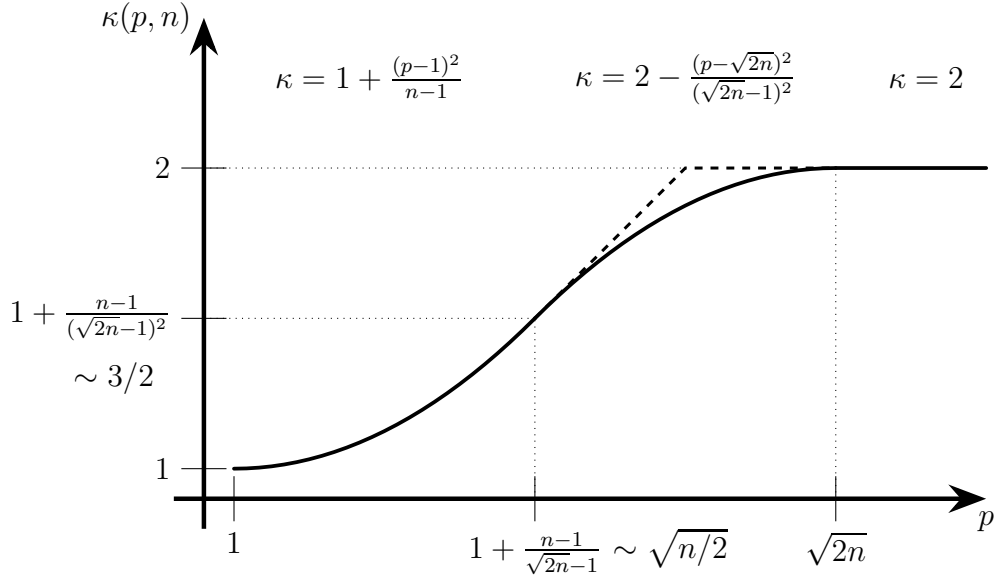
\begin{figure}[h!]
%\documentclass[tikz,border=1mm]{standalone}
%\usetikzlibrary{calc, decorations.pathreplacing, arrows.meta}
%\usepackage{amsmath}
%\begin{document}

% My favorite arrow style
\tikzstyle{mmArrow}=[-{Stealth[length=#1]}]
% style for drowing the plot line
\tikzstyle{PlotLine} = [line width=0.5mm, line cap=round]
% styles for drawing ticks on axes, with number written next to it
\tikzset{hTick/.pic={
	\draw (-0.15,0) node[left] {#1} -- ++(0.3,0);
}}
\tikzset{vTick/.pic={
	\draw (0,-0.15) node[below] {#1} -- ++(0,0.3);
}}
% over brace
\tikzstyle{bzBraceO} = [draw,
	decoration={brace, raise=0.0cm, amplitude=3mm},
	decorate,
	every node/.style={anchor=south, yshift= 0.2cm}]

\begin{tikzpicture}[x=2cm, y=2cm]
	% The origin is shifted by epsD to makes certain other things easier
	\def\epsD{0.2}
	% axes
	\draw[mmArrow=5mm, line width=0.6mm, shift={(0,-\epsD)}] (-0.4,0) -- (5,0) node[below] {$p$};
	\draw[mmArrow=5mm, line width=0.6mm, shift={(-\epsD,0)}] (0,-0.4) -- (0,3) node[left] {$\kappa(p,n)$};
	% axes labels
	\pic at (-\epsD,0) {hTick=$1$};
	\pic at (-\epsD,1) {hTick=$1+\tfrac{n-1}{(\sqrt{2n}-1)^2}$};
	\node at (-0.8,0.6) {$\sim 3/2$};
	\pic at (-\epsD,2) {hTick=$2$};
	\pic at (0,-\epsD) {vTick=$1$};
	\pic at (2,-\epsD) {vTick=$1+\tfrac{n-1}{\sqrt{2n}-1}$};
	\node at (2.9,-0.53) {$\sim \sqrt{n/2}$};
	\pic at (4,-\epsD) {vTick=$\sqrt{2n}$};

	% plots with labels
	\draw[line width=0.4mm, dashed] (2,1) -- (3,2) -- (4,2);	% dashed line showing the graph for the case k=1
  	\draw[PlotLine, domain=0:2]
  		plot (\x,{\x*\x/4});
  	\draw[PlotLine, domain=2:4]
  		plot (\x,{2-(\x-4)*(\x-4)/4});
  	\draw[PlotLine]
  		(4,2) -- (5,2);

  	% values in the three regions
  	\node[scale=1.0] at (0.9,2.6) {$\kappa = 1+\tfrac{(p-1)^2}{n-1}$};
  	\node[scale=1.0] at (3.0,2.6) {$\kappa = 2-\tfrac{(p-\sqrt{2n})^2}{(\sqrt{2n}-1)^2}$};
  	\node[scale=1.0] at (4.6,2.6) {$\kappa = 2$};

  	% help lines from the x axis (or rather p axis) to the plot
	% and also from the y axis (or rather kappa axis) to the plot
  	\draw[dotted]
  		(2,-\epsD) -- (2,1)
  		(4,-\epsD) -- (4,2)
		(-\epsD,1) -- (2,1)
		(-\epsD,2) -- (3,2);
\end{tikzpicture}
%\end{document} 
\caption{A graph showing the dependence of $\kappa(p,n)$ on $p$ (assuming $d \ge 2$). The dashed curve illustrates the best constant $\kappa(p,n)$ in the scalar case $d = 1$ \cite{ChaCheWei16}. Note that the shape presented above corresponds to the asymptotic case $n \to \infty$, where the first two regimes have the same length. For small $n$, the middle regime can be very little (its length is $\approx 0.07$ for $n = 3$) or even absent (for $n = 2$, as found in \cite{MM}). In general, $\kappa(p,n)$ is a piecewise quadratic function of class $C^1$.}
\end{figure}

In \Cref{ch:application-of-Kato} we show how our Kato type inequality can be used to study the regularity of $p$-harmonic maps. As already mentioned in this introduction, the reasoning is classical and follows by combining a stability inequality with a Bochner formula, as in \cite{SU3}. For suitably chosen $p$, $n$ and $d$, this leads to full or partial regularity of minimizing $p$-harmonic maps $B^n \to \S^d$. In particular, if $n$ and $d$ are chosen so that harmonic maps from $\B^n$ into $\S^d$ are known to be regular (see \cite{MM} for a summary of all such cases), then we provide an interval of $p \in [2,p_0(n,d)]$ for which $p$-harmonic maps are also regular. To give a concrete example, we study the case of $p$-harmonic maps $B^4 \to S^4$ in detail, obtaining the interval $p \in [2,p_0]$ with $p_0 = \frac{5+\sqrt{13}}{3} \approx 2.86$.

\Cref{ch:regularity-B3-S3} provides further improvement in the special case of $p$-harmonic maps $B^3 \to \S^3$, which was the main focus in \cite{MM}. There regularity was established for all $2 \le p \le 2.366$ and $p \ge 2.962$. Here, the main innovation lies in a new choice of test functions for both the stability and Bochner inequalities: instead of $|\nabla u|$ and $|\nabla u|^{p-2}$ (respectively), which is the natural choice, one can take $|\nabla u|^{1+\gamma}$ and $|\nabla u|^{p-2+2\gamma}$ with an additional parameter~$\gamma$. By formal calculation, one can check that choosing $\gamma$ strictly below zero leads to a better constant. However, an optimal choice $\gamma \approx -0.6$ is not admissible, as the exponent $p-2+2\gamma$ becomes negative in this case. For this reason, we perform a non-trivial regularization procedure, and finally obtain regularity of minimizing $p$-harmonic maps $B^3 \to \S^3$ for all $2 \le p \le 2.642$ as a result. 

{\bf Acknowledgments.} The project is co-financed by:
\begin{itemize}
\item the Deutsche Forschungsgemeinschaft, project no. 441380936 (AG);
 \item  the Polish National Agency for Academic Exchange within Polish Returns Programme - BPN/PPO/2021/1/00019/U/00001 (KM);
 \item the National Science Centre, Poland grant 2022/01/1/ST1/00021 (KM);
 \item the National Science Centre, Poland grant 2019/35/B/ST1/02030 (MM).
\end{itemize}

\section{Kato inequality}
\label{ch:Kato}
We may assume that $\mathcal N$ is isometrically embedded into some Euclidean space $\R^L$. We recall the Euler--Lagrange equations for $p$-harmonic maps $u\in W^{1,p}(\cM,\cN)$
\begin{equation}\label{eq:EL}
- \div_\cM(|\nabla u|^{p-2} \nabla u) = |\nabla u|^{p-2} A(\nabla u, \nabla u),
\end{equation}
where $A$ is the second fundamental form of the isometric embedding $\mathcal N\subset \R^L$.

\begin{lemma}[an elementary fact from Linear Algebra]\label{le:elementary}
  For any linear map $A \in L(\R^n, \R^d)$, 
  there is an
  orthonormal basis $(b_1,\ldots,b_n)$ of $\R^n$ for which
  $(Ab_1,\ldots,Ab_n)$ are pairwise orthogonal in $\R^d$.
\end{lemma}

\begin{proof} Let $(e_1,\ldots,e_n)$ be the canonical basis of $\R^n$. We need a $B\in SO(n)$ such that $\la ABe_i,ABe_j\ra=\lambda_i\delta_{ij}$ for all $i,j\in\{1,\ldots,n\}$. This is equivalent to $\la B^{T}A^{T}ABe_i,e_j \ra=\lambda_i\delta_{ij}$, which means all we have to do is diagonalize $A^{T}A\in L(\R^n,\R^n)$. But, since $A^{T}A$ is symmetric, such $B\in SO(n)$ exists. Now $b_i\coloneqq Be_i$ proves the assertion.
\end{proof}

Having this, we prove
\begin{theorem}[Kato's inequality for $p$-harmonic maps]\label{thm:Kato}
  Let $\mathcal M$, $\mathcal N$ be smooth Riemannian manifolds of dimensions $n\ge2$ and $d\ge2$. Let $u\colon \m\to \n$ be a $p$-harmonic map, $p>1$. Then, at any point $q\in \m$
  where $\nabla u(q)\ne0$, we have
  \begin{equation}\label{eq:Katoinequality}
    \kappa(p,n)\,|\nabla|\nabla u||^2\le|\nabla^2u|^2,
  \end{equation}
  where
  \begin{equation}\label{eq:kappaformula}
    \arraycolsep=1.4pt\def\arraystretch{2.2}
    \kappa(p,n) =\left\{ \begin{array}{lll}
                     1+\frac{(p-1)^2}{n-1} &\quad  \text{for }&1\le p\le1 + \frac{n-1}{\sqrt{2n}-1}\\
                        2-\frac{(p-\sqrt{2n})^2}{(\sqrt{n}-\sqrt{2})^2} &\quad \text{for }&1 + \frac{n-1}{\sqrt{2n}-1}\le p \le \sqrt{2n}\\
                     2 &\quad \text{for }& p\ge \sqrt{2n}
                    \end{array}
\right.
  \end{equation}
  is the sharp constant in \eqref{eq:Katoinequality}, in the sense that there is map that solves the $p$-harmonic equation at $q$ and achieves equality at that point.
\end{theorem}

\begin{remark}\label{re:Katoindim2}
Note that for $n=2$ in the formula \eqref{eq:kappaformula} we have the equality
\[
 1+ \frac{n-1}{\sqrt{2n}-1} = 2 = \sqrt{2n}
\]
and hence the formula corresponds to the one in \cite{MM}, i.e., $\kappa(p,2) = \min\left\{1+(p-1)^2,2\right\}$.
\end{remark}

\begin{proof}[Proof of \Cref{thm:Kato}]
From now, we fix a point $q \in \cM$ with $\nabla u(q) \neq 0$ and use normal coordinates around $q \in \cM$ and $u(q) \in \cN$. Then, the $p$-harmonic equation at this particular point takes the Euclidean form.
%We use the $p$-harmonic equation in one point $q\in \m$ where the coordinates are ``Euclidean''. 
That is, for $a\in\{1,\ldots,d\}$,
\begin{equation}\label{eq:plaplace}
  \sum_{i=1}^n u_{ii}^a
  +\frac{p-2}{|\nabla u|^2}\sum_{i,j=1}^n\sum_{b=1}^d u_i^a u_j^b u_{ij}^b=0.  
\end{equation}
By \Cref{le:elementary}, we may choose orthonormal bases of $T_q \m$ and $T_{u(q)} \cN$ in which $\nabla u$ has only diagonal components. In other words, a linear change of coordinates ensures that 
\begin{equation}\label{eq:coordinates-choice}
u_i^a = \alpha_i \delta_i^a
  \quad 
\text{for all } i \in \{ 1, \ldots, n \}, \ a\in\{1,\ldots,d\}. 
\end{equation}
From now on, we will additionally assume $|\nabla u(q)|^2 = 1$, which translates to the condition $\sum_{a=1}^n \alpha_a^2 = 1$. Otherwise, we would just work with coefficients $\tilde{u}_i^a = u_i^a / |\nabla u|$ instead of $u_i^a$. Moreover, for notational convenience we shall adopt the convention that $u_{ij}^a = 0$ whenever $a > d$. 
Thus, the $p$-harmonic system \eqref{eq:plaplace} takes the form
\begin{equation}\label{eq:plaplacecoordinatechoice}
 \Delta u^a + (p-2) \alpha_a \sum_{b=1}^n \alpha_b u_{ab}^b = 0 
 \quad \text{for each } a \in \{1,\ldots, d\}. 
\end{equation}

Using \eqref{eq:coordinates-choice}, we introduce the quantity 
\begin{equation}
I_i := |\pl_i |\nabla u||^2 
  = \brac{\sum_{j=1}^n\sum_{b=1}^n u_j^b u_{ij}^b}^2 
  = \brac{\sum_{b=1}^n \alpha_b u_{ib}^b}^2 
  \quad \text{for } i \in \{ 1,\ldots,n \}
\end{equation}
and aim to show the estimate 
\begin{equation}\label{eq:Ii}
\begin{split}
 \kappa(p,n)I_i 
 \le \sum_{j=1}^n (u_{jj}^i)^2 + 2\sum_{b\neq i}(u_{ib}^b)^2.
 \end{split}
\end{equation}

Once we establish \eqref{eq:Ii}, the Theorem follows, since then 
\begin{equation}
\label{eq:hessiandivided}
\begin{split}
 |\nabla^2 u|^2 & = \sum_{i,j=1}^n\sum_{a=1}^d (u_{ij}^a)^2 
% = \sum_{i=d+1}^n\sum_{b=1}^d\brac{(u_{ib}^b)^2+(u_{bi}^b)^2} + 
 \ge \sum_{i=1}^n \brac{\sum_{j=1}^n (u_{jj}^i)^2 + 2\sum_{b\neq i} (u_{ib}^b)^2} \\
 & \ge \sum_{i=1}^n \kappa(p,n) I_i 
 = \kappa(p,n) |\nabla |\nabla u||^2. 
\end{split}
\end{equation}
Note the first inequality becomes an equality if and only if $u_{ij}^a = 0$ for all triples $i,j,a$ of pairwise distinct indices.
Thus, it suffices to show \eqref{eq:Ii}. For simplicity of notation we will focus on the case when $i=1$ (the remaining cases are proved in the same way). 

First, we rewrite the equation \eqref{eq:plaplacecoordinatechoice} with $a = 1$ in the following way: 
\begin{equation}\label{eq:p-harmeqwith1}
 \sum_{j=2}^n u_{jj}^1 = -u_{11}^1 - (p-2)\alpha_1\sum_{b=1}^n \alpha_b u_{1b}^b.
\end{equation}
Combining this with the elementary inequality between the quadratic and arithmetic means, we obtain
\begin{equation}\label{eq:1overn-1}
 \sum_{j=2}^n(u_{jj}^1)^2 \ge \frac{1}{n-1} \brac{\sum_{j=2}^n u_{jj}^1}^2 = \frac{1}{n-1} \brac{u_{11}^1 + (p-2)\alpha_1\sum_{b=1}^d \alpha_b u_{1b}^b}^2.
\end{equation}
Note that equality in \eqref{eq:1overn-1} occurs when all terms $u_{jj}^1$ coincide for $j \in \{2,\ldots,n\}$. 

Applying \eqref{eq:1overn-1}, we see that \eqref{eq:Ii} reduces to showing 
\begin{equation}\label{eq:Iinew}
\kappa(p,n)\brac{\sum_{b=1}^n \alpha_b u_{1b}^b}^2 \le  (u_{11}^1)^2 + \frac{1}{n-1} \brac{u_{11}^1 + (p-2)\alpha_1\sum_{b=1}^n \alpha_b u_{1b}^b}^2 + 2\sum_{b=2}^n (u_{1b}^b)^2.
\end{equation}
As a consequence, our problem now simplifies to the elementary inequality above, which we shall show is valid for any $\alpha_i$ ($i \in \{ 1,\ldots,n \}$) satisfying $\alpha_1^2+\ldots+\alpha_n^2 = 1$, and any $u_{1b}^b$ ($b \in \{ 1,\ldots,n \}$).

The next step can be seen as a further reduction to the case $d = 2$. Let us rewrite the above, extracting all Hessian terms other than $u_{11}^1$: 
\begin{equation}\label{eq:Iinewer}
\begin{split}
\kappa(p,n)&\brac{\alpha_1u_{11}^1 + \sum_{b=2}^d \alpha_b u_{1b}^b}^2\\
&\le  (u_{11}^1)^2 + \frac{1}{n-1} \brac{u_{11}^1 + (p-2)\alpha_1\brac{\alpha_1u_{11}^1+\sum_{b=2}^d \alpha_b u_{1b}^b}}^2 + 2\sum_{b=2}^d(u_{1b}^b)^2.
\end{split}
\end{equation}
Let us denote $\beta = |(\alpha_2,\ldots,\alpha_d)|$ and choose $y$ satisfying $\beta y = \sum_{b=2}^d \alpha_b u_{1b}^b$ (if $\beta = 0$, choose $y = 0$). By the Cauchy--Schwarz inequality, $(\beta y)^2 \le \beta^2 \sum_{b=2}^d(u_{1b}^b)^2$, which we can rewrite as $y^2 \le \sum_{b=2}^d(u_{1b}^b)^2$. Thus, after denoting $x := u_{11}^1$ and $\alpha := \alpha_1$, our inequality reduces to 
\begin{equation}\label{eq:Iinewererer}
\kappa(p,n)\brac{\alpha x+\beta y}^2 \le  x^2 + \frac{1}{n-1} \brac{x + (p-2)\alpha\brac{\alpha x+\beta y}}^2 + 2y^2.
\end{equation}
which holds, as we shall see, for all $x,y,\alpha,\beta$ satisfying $\alpha^2 + \beta^2 = 1$. 
The appearance of $\alpha x + \beta y$ suggests a change of variables (a rotation), which we now perform: 
\begin{equation}
\label{eq:change-xy-zw}
\begin{cases}
     z = \alpha x+\beta y, \\
     w = \beta x-\alpha y,
\end{cases}
\quad \text{ or } \quad 
\begin{cases}
x = \alpha z+\beta w, \\
y = \beta z-\alpha w.
\end{cases}
\end{equation}
Inequality \eqref{eq:Iinewererer} now becomes
\begin{equation}\label{eq:Iinewerererer}
\begin{split}
0 &\le  z^2 \brac{ -\kappa(p,n) + 2 + \alpha^2 \left( \frac{(p-1)^2}{n-1} - 1 \right) } 
+ w^2\brac{2-\frac{n-2}{n-1}\beta^2} - zw\frac{2\alpha\beta(n-p)}{n-1}.
\end{split}
\end{equation}
Now \eqref{eq:Iinewerererer} is a quadratic form in $z$ and $w$, with a positive coefficient of $w^2$. Hence the inequality is univerally true if only its discriminant $\Delta_{zw}$ is nonpositive, that is
\begin{equation}\label{eq:determinantinzw}
 \begin{split}
 0 \ge \Delta_{zw}
 &= \brac{\frac{2\alpha\beta(n-p)}{n-1}}^2 
 - 4 \brac{ -\kappa(p,n) + 2 + \alpha^2 \left( \frac{(p-1)^2}{n-1} - 1 \right) }
 \brac{2-\frac{n-2}{n-1}\beta^2}\\
 &= -\frac{4}{n-1}\brac{-\kappa(p,n)((n-2)\alpha^2+n)+\alpha^4(p-2)^2+\alpha^2(p^2-4)+2n}.
\end{split}
 \end{equation}
Note that we have used the equality $\alpha^2+\beta^2 = 1$ to eliminate $\beta$. Now \eqref{eq:determinantinzw} is equivalent to 
\begin{equation}\label{eq:determinantasinequality}
 \kappa(p,n) \le \frac{\alpha^4(p-2)^2+\alpha^2(p^2-4)+2n}{(n-2)\alpha^2+n}
\end{equation}
for all $\alpha\in[-1,1]$. Substituting $a := \alpha^2 \in [0,1]$, we obtain
\begin{equation}\label{eq:kappaitisjustminimize}
 \kappa(p,n) = \min_{a\in[0,1]}\frac{a^2(p-2)^2+a(p^2-4)+2n}{(n-2)a+n}\eqqcolon \min_{a\in[0,1]}f(a).
\end{equation}
Minimization of $f(a)$ is straightforward, but its result depends --- quite unsurprisingly --- on $n$ and $p$. Let us compute:
\begin{align}
  f'(a) &= \frac{a^2(n-2)(p-2)^2+2an(p-2)^2+n(p^2-2n)}{(a(n-2)+n)^2}.\label{eq:derivativeoff}
\end{align}
Note that the numerator of $f'(a)$ is a quadratic polynomial with nonnegative discriminant $\Delta_a = 8n(p-2)^2(n-p)^2$. Thus, $f'(a)$ is nonpositive exactly on an interval $a \in [a_-, a_+]$ satisfying $a_+ + a_- = -2n(p-2)^2$ and $a_+ \cdot a_- = n(p^2-2n)$. Naturally, the behavior of $f$ on the interval $[0,1]$ is different in the following three cases:

\textsc{Case 1:} $f'(0) \ge 0$. This happens when $p \ge \sqrt{2n}$. In that case, we see that $a_-,a_+$ are both nonpositive, so $f$ is increasing on $[0,1]$. In conclusion: 
\[
\kappa(p,n) = f(0) = 2 
\qquad \text{for }
p \ge \sqrt{2n}. 
\]

\textsc{Case 2:} $f'(1) \le 0$. If we restrict to $p \ge 1$, this happens when $p \le 1 + \frac{n-1}{\sqrt{2n}-1}$. In that case, we see $a_- \le 0 \le 1 \le a_+$, so $f$ is decreasing on $[0,1]$, and in consequence: 
\[
\kappa(p,n) = f(1) = 1+\frac{(p-1)^2}{n-1} 
\qquad \text{for }
1 \le p \le 1 + \frac{n-1}{\sqrt{2n}-1}. 
\]

\textsc{Case 3:} $f'(0) < 0 < f'(1)$. By elimination, this is true in the regime $1 + \frac{n-1}{\sqrt{2n}-1} < p < \sqrt{2n}$. In that case, we have $a_+ \in (0,1)$ and this is where the minimum is achieved. We have already computed the discriminant -- since in this regime $p \in [2,n]$, we have 
\begin{align*}
\sqrt{\Delta_a} & = 2\sqrt{2n}(p-2)(n-p), \\
a_+ & = \frac{-n(p-2) + \sqrt{2n}(n-p)}{(n-2)(p-2)}, \\
\kappa(p,n) = f(a_+) & = 2-\frac{(p-\sqrt{2n})^2}{(\sqrt{n}-\sqrt{2})^2}
\qquad \text{in the middle regime.}
\end{align*}
This analysis completes the proof of the Kato inequality.

\medskip

Finally, we complete the proof of \Cref{thm:Kato} by showing that the constant $\kappa(p,n)$ defined in \eqref{eq:kappaformula} is sharp. For simplicity, we shall do so for a map $u \colon \R^n \to \R^2$. One can then obtain an $\R^d$-valued map by the inclusion $\R^2 \cong \R^2 \times \{ 0 \} \subseteq \R^d$, and a map between manifolds via normal coordinates (as in the beginning of the proof). To find the values of the derivatives $u_i^a$ and $u_{ij}^a$ (at $0$), we work backwards so that each step is optimal: 
\label{SHARPNESS}
\begin{itemize}
\item Take $a \in [0,1]$ for which the minimum $\min_{a \in [0,1]} f(a)$ in \eqref{eq:kappaitisjustminimize} is achieved, and set $\alpha := \sqrt{a}$, $\beta := \sqrt{1-\alpha^2}$. In the case $d = 2$, we then simply put 
\[
u_1^1 := \alpha_1 = \alpha, \qquad 
u_2^2 := \alpha_2 = \beta,
\]
and all other first order derivatives equal zero. 
\item By this choice, we have $\Delta_{zw} = 0$ in \eqref{eq:determinantinzw}. This means that we can pick a nonzero pair $(z,w)$ for which equality in \eqref{eq:Iinewerererer} is achieved. Reversing the change of variables in \eqref{eq:change-xy-zw}, we obtain a nonzero pair $(x,y)$ achieving equality in \eqref{eq:Iinewererer}. 
\item Choose $u_{11}^1 = x$, $u_{12}^2 = y$. For $d = 2$, the Cauchy--Schwarz inequality is trivially an equality. In consequence, \eqref{eq:Iinewer} and \eqref{eq:Iinew} are also equalities. 
\item In order to satisfy the $p$-harmonic equation \eqref{eq:p-harmeqwith1}, we choose 
\[
u_{22}^1,\ldots,u_{nn}^1 := - \frac{1}{n-1} \left( u_{11}^1 + (p-2) \alpha_1 \sum_{b=1}^n \alpha_b u_{1b}^b \right).
\]
That way, the quadratic--arithmetic mean inequality between these terms \eqref{eq:1overn-1} is an equality, which means that \eqref{eq:Ii} is, too: 
\[
\kappa(p,n) |\pl_1 |\nabla u||^2 
= \sum_{j=1}^n (u_{jj}^i)^2 + 2 \sum_{b\neq i}(u_{ib}^b)^2.
\]
\item We let all other Hessian terms be zero, which makes the $p$-harmonic system \eqref{eq:coordinates-choice}, \eqref{eq:plaplacecoordinatechoice} satisfied on all coordinates. Also, the equality above is then actually the equality between $\kappa(p,n) |\nabla |\nabla u||^2$ and $|\nabla^2 u|^2$, as all other terms on both sides vanish. 
\end{itemize}
\end{proof}

\begin{remark}
 The $p$-harmonic map $u\colon \R^n \to \R$, $u(x)\coloneqq |x|^{\frac{p-n}{p-1}}$, realizes the bound $\kappa (p,n) \le 1 + \frac{(p-1)^2}{n-1}$. It can be concatenated with any geodesic to give $p$-harmonic maps from $\R^n$ to $\mathcal N$.
 
 The mapping $v\colon \R^n \to \S^d$, $v(x)\coloneqq \frac{(x_1,\ldots,x_R,0,\ldots,0)}{\sqrt{x_1^2+\ldots x_R^2}}$ is $p$-harmonic for all $p>1$ and all $R\le \min\{d+1,n\}$, and it realizes the bound $\kappa(p,n)\le 2$ in all relevant dimensions.
 
 These examples show that the first and the third case of \Cref{thm:Kato} are optimal, not only in the pointwise sense. We would expect the same for the second case.
\end{remark}

\section{Application of Kato's inequality to regularity of \texorpdfstring{$p$}{p}-harmonic maps}
\label{ch:application-of-Kato}

As an application of the vectorial Kato inequality, we can obtain regularity results for $p$-harmonic mappings into spheres.
\begin{theorem}\label{th:regularityfromKato}
Assume that $n \ge 3$ and $2 \le p < d$, with $n$, $d$ both integers. If $n \le n_0(p,d)$, then every minimizing $p$-harmonic map $u\in W^{1,p}(B^n,\S^d)$ is regular in the interior of $B^n$. If $n = n_0(p,d)+1$, then every minimizing $p$-harmonic map $u\in W^{1,p}(B^n,\S^d)$ has at most isolated singularities. In general the singular set of a minimizing $p$-harmonic map $u\in W^{1,p}(B^n,\S^d)$ is a closed set of Hausdorff dimension at most $n-n_0(p,d)-1$.
 
For fixed $p$ and $d$, the number $n_0(p,d)$ is defined as the largest integer $n$ satisfying the inequalities 
\begin{align}
\label{eq:Kato=>regularity}
%\label{eq:Kato=>regularity-imp}
\kappa(p,n-1) & \ge \frac{(d+p-2)(n-2)}{(d-p)(n-1)} - (p-2), \\
\label{eq:Kato=>regularity-non-imp}
d(n-p)^2 & < 4(d-p)(n-1).
\end{align}
\end{theorem}

\begin{remark}
\label{rem:n0-properties}
Observe the following properties: 
\begin{enumerate}[label=(\alph*)]
\item \label{rem:ita}
If conditions \eqref{eq:Kato=>regularity}, \eqref{eq:Kato=>regularity-non-imp} are satisfied for $n = n_0$ (with $p$, $d$ fixed), then they are also satisfied for any $p < n \le n_0$. 
\item \label{rem:itb}
If conditions \eqref{eq:Kato=>regularity}, \eqref{eq:Kato=>regularity-non-imp} are satisfied for some $d$ (with $p$, $n_0$ fixed), then they are also satisfied for any larger $d$. 
\end{enumerate}
\end{remark} 

\begin{proof}[Proof of Remark \ref{rem:n0-properties}]
Item \ref{rem:ita}. It follows from the formula \eqref{eq:kappaformula} for $\kappa(p,n)$ that it is a~non-increasing function of $n$. On the other hand, $\frac{n-2}{n-1}$ is increasing for $n > 1$, so \eqref{eq:Kato=>regularity} is true for all $1 < n \le n_0$. As for \eqref{eq:Kato=>regularity-non-imp}, it is trivially true for $n = p$, so it also holds for all $p \le n \le n_0$ by convexity of the left-hand side. 

Item \ref{rem:itb}. We simply note that $\kappa(p,n-1)$ does not depend on $d$, whereas both functions $\frac{d+p-2}{d-p}$ and $\frac{d}{d-p}$ are decreasing in $d$. 
\end{proof}

In the proof of \Cref{th:regularityfromKato}  we will follow the general strategy of Schoen and Uhlenbeck \cite{SU3}, let us begin with the defition of a tangent map.
\begin{definition}\label{de:tangentmap}
A $p$-harmonic map $u\colon \R^n \to \n$ is called \emph{tangent} if it is radially symmetric, i.e., $u(x) = \tilde u(\frac{x}{|x|})$ for a mapping $\tilde u\colon \S^{n-1}\to\n$. Since $u$ is uniquely determined by $\tilde u$ it will be useful for us to treat tangent maps as maps defined on the sphere without changing the notation\footnote{Of course, since the map $u$ is $p$-harmonic, the map $\tilde u\colon \S^{n-1}\to\n$ satisfies the Euler--Lagrange equation for $p$-harmonic maps from $\S^{n-1}$ to $\n$.}.
If in addition, $u$ minimizes the $E_p$ energy on compact subsets of $\R^n$, we call it a \emph{minimizing tangent map}.
\end{definition}

We will use the following two inequalities for tangent maps.

\begin{lemma}\label{le:stability&bochner}
Assume that $p \ge 2$ and $u\colon \S^{n-1} \to \S^d$ is a minimizing $p$-harmonic tangent map. Then the following estimates hold:
\begin{itemize}
 \item[(a)] Stability inequality:
 \begin{equation}\label{eq:stabilty}
\begin{split}
\int_{\S^{n-1}} &|\nabla u|^{p-4} \left( d |\nabla u|^2 |\nabla |\nabla u||^2 + (p-2) |\la \nabla u, \nabla |\nabla u| \ra|^2 \right)\\
&\ge (d-p) \int_{\S^{n-1}} |\nabla u|^{p+2} -  \frac{d(n-p)^2}{4} \int_{\S^{n-1}} |\nabla u|^p; 
\end{split}
\end{equation}
\item[(b)] Bochner inequality:
\begin{equation}\label{eq:bochner}
\begin{split}
\int_{\S^{n-1}} &|\nabla u|^{p-2} (|\nabla^2 u|^2 + (p-2) |\nabla |\nabla u||^2)\\
&\le \frac{n-2}{n-1} \int_{\S^{n-1}} |\nabla u|^{p+2} - (n-2) \int_{\S^{n-1}} |\nabla u|^p. 
\end{split}
\end{equation}
\end{itemize}
\end{lemma}

The proof of the Stability inequality \eqref{eq:stabilty} essentially follows from \cite[Lemmata 2 and 3]{Nakauchi96}. Similarly, the proof of Bochner's inequality \eqref{eq:bochner} also relies on Nakauchi \cite[Lemma 1]{Nakauchi96} (see also \cite[Lemma 2]{Nakauchi01}), with a minor modification: omitting the step $|\nabla|\nabla u|| \le |\nabla^2 u|$. For the reader's conveniencee, we provide the proofs with the necessary modifications in \Cref{ch:appendix}.

\begin{proof}[Proof of \Cref{th:regularityfromKato}]
By \cite[Theorem 4.5]{HLp} it suffices to prove that if $n \le n_0(p,d)$, then every minimizing tangent map $\bar{u} \colon \S^{n-1} \to \S^d$ must be constant. Note that the cases $n \le p$ already follow from the general result \cite[Corollary 2.6]{HLp}, valid for any target manifold. Thus, we only consider $p < n \le n_0(p,d)$ in the following. Recall that both conditions \eqref{eq:Kato=>regularity} and \eqref{eq:Kato=>regularity-non-imp} are still satisfied for $n$ by \Cref{rem:n0-properties}. 

We apply \Cref{le:stability&bochner} to $\bar{u}$. Combining the Bochner inequality \eqref{eq:bochner} with the Kato inequality \Cref{thm:Kato} $|\nabla^2 \bar u|^2 \ge \kappa |\nabla |\nabla \bar u||^2$, we obtain
\begin{equation}\label{eq:bochner2}
\begin{split}
(\kappa + p-2)&\int_{\S^{n-1}} |\nabla \bar u|^{p-2} |\nabla |\nabla \bar u||^2\\
&\le \frac{n-2}{n-1} \int_{\S^{n-1}} |\nabla \bar u|^{p+2} - (n-2) \int_{\S^{n-1}} |\nabla \bar u|^p. 
\end{split}
\end{equation}
Note that the largest possible $\kappa$ is $\kappa(p,n-1)$ given by \eqref{eq:Katoinequality} in \Cref{thm:Kato}, but we may work with a smaller constant as well. Next, applying the Cauchy--Schwarz inequality to the inner product appearing in the Stability inequality \eqref{eq:stabilty}, we have
\begin{equation}\label{eq:stabilty2}
\begin{split}
\int_{\S^{n-1}} |\nabla \bar u|^{p-2} |\nabla |\nabla \bar u||^2 
\ge \frac{d-p}{d+p-2} \int_{\S^{n-1}} |\nabla \bar u|^{p+2} -  \frac{d(n-p)^2}{4(d+p-2)} \int_{\S^{n-1}} |\nabla \bar u|^p. 
\end{split}
\end{equation}
Combining the last two inequalities, we can eliminate the second order terms, which leads to 
\begin{equation}
\begin{split}
 &\brac{\frac{(d-p)(\kappa+p-2)}{d+p-2}-\frac{n-2}{n-1}} \int_{\S^{n-1}} |\nabla \bar u|^{p+2} \\
 &\quad+ \brac{n-2 -  \frac{d(n-p)^2(\kappa+p-2)}{4(d+p-2)}}\int_{\S^{n-1}} |\nabla \bar u|^p \le 0. 
 \end{split}
\end{equation}
Now, we will show that both terms in the brackets are nonnegative, and the second one is nonzero. In consequence, we can conclude that $\int_{\S^{n-1}} |\nabla \bar u|^p = 0$, which means that $\bar u$ is constant, as required. 

Thus, we proceed to check the sign of these two terms. To this end, let us rewrite our claim in terms of $\kappa$: 
\begin{equation}\label{eq:joint-reg-condition}
E_1 := \
 \frac{(n-2)(d+p-2)}{(n-1)(d-p)} \le \kappa + p-2 \le \frac{4(n-2)(d+p-2)}{d(n-p)^2}
\ =: E_2.
\end{equation}
As we already remarked, the bounds \eqref{eq:Kato=>regularity} and \eqref{eq:Kato=>regularity-non-imp} hold for our choice of $p$, $d$ and $n$. In the present notation, these two bounds are equivalent to $\kappa(p,n-1) + p-2 \ge E_1$ and $E_1 \le E_2$, respectively. Thus, we can choose 
\[
\kappa := \min \left( \kappa(p,n-1), E_1 - (p-2) \right),
\]
so that our claim holds. This completes the proof of \Cref{th:regularityfromKato}. 
\end{proof}

Although \Cref{th:regularityfromKato} gives an explicit condition \eqref{eq:Kato=>regularity}, it is still difficult to see for which $p$ (depending on $n,d$) we actually obtain regularity. To illustrate how such analysis can be performed, we state a separate corollary for $n = d = 4$.
\begin{corollary}
%For $n=d=4$ we obtain that 
Minimizing $p$-harmonic mappings $u \in W^{1,p}(B^4,\S^4)$ are regular whenever 
\[
 p \in [2,p_0], \quad \text{ where } \quad 
 p_0 = \frac{5+\sqrt{13}}{3} \approx 2.8685.
\]
\end{corollary}

\begin{proof}
By \Cref{th:regularityfromKato}, it is enough to check that the conditions \eqref{eq:Kato=>regularity} and \eqref{eq:Kato=>regularity-non-imp} are satisfied for $n = d = 4$ and $p \in [2,p_0]$. In our case, these conditions read 
\[
\kappa(p,3) \ge \frac{2(p+2)}{3(4-p)} - (p-2), 
\qquad 
4 \cdot (4-p)^2 < 4 \cdot (4-p) \cdot 3,
\]
respectively. The latter condition reduces to $p > 1$, so we now focus on the former. We employ the Kato constant from \Cref{thm:Kato}:
\[
\kappa(p,3) = \left\{ \begin{array}{lll}
                     1+\frac{(p-1)^2}{2} &\quad  \text{for }&1\le p\le \frac{7+2\sqrt{6}}{5}\\
                        2-\frac{(p-\sqrt{6})^2}{(\sqrt{3}-\sqrt{2})^2} &\quad \text{for }&\frac{7+2\sqrt{6}}{5}\le p \le \sqrt{6}\\
                     2 &\quad \text{for }& p\ge \sqrt{6}.
                    \end{array}
\right.
\]
One can easily graph the functions $\kappa(p,3)$ and $\frac{2(p+2)}{3(4-p)} - (p-2)$ to see when the desired inequality occurs.
\begin{figure}[h!]
%\documentclass[tikz,border=1mm]{standalone}
%\usetikzlibrary{calc, decorations.pathreplacing, arrows.meta}
%\usepackage{amsmath}
%\begin{document}

% My favorite arrow style
\tikzstyle{mmArrow}=[-{Stealth[length=#1]}]
% style for drowing the plot line
\tikzstyle{PlotLine} = [line width=0.3mm, line cap=round]
% styles for drawing ticks on axes, with number written next to it
\tikzset{hTick/.pic={
	\draw (-0.15,0) node[left] {#1} -- (0,0);
}}
\tikzset{vTick/.pic={
	\draw (0,-0.15) node[below] {#1} -- (0,0);
}}

\begin{tikzpicture}[scale=3.5]
	% The origin is shifted to make certain other things easier
	\def\Xshift{1.8}
	\def\Yshift{1.2}
	\tikzstyle{OrShift}=[shift={(\Xshift,\Yshift)}]
	% help lines
	\tikzstyle{Help} = [thin, gray, dashed]
	\draw[Help]
		(2,\Yshift) -- (2,1.5)
		(2.5,\Yshift) -- (2.5,2)
		(\Xshift,1.5) -- (2.5,1.5)
		(\Xshift,2) -- (2.45,2);
	% \frac{7+2\sqrt{6}}{5} = 2.38
	% \sqrt{6} = 2.45
	% 1/(\sqrt{3}-\sqrt{2}) = 9.899
	% the actual plot lines
	\draw[PlotLine]
  		plot[domain=2:2.38] (\x, {1 + (\x-1)*(\x-1)/2})
  		plot[domain=2.38:2.45] (\x, {2 - 9.899*(\x-2.45)*(\x-2.45)})
  		-- (3,2);
	\draw[PlotLine]
		plot[domain=2:3] (\x, {0.667*(\x+2)/(4-\x) - (\x-2)});
	% explanation for plots 
	\node[anchor=south] at (2.5,2.05) {$\kappa(p,3)$};
	\node[anchor=west] at (2.6,1.5) {$\tfrac{2(p+2)}{3(4-p)} - (p-2)$};
	% axes
	\draw[mmArrow=2.5mm, line width=0.3mm, OrShift] (-0.1,0) -- ++(1.5,0);
	\draw[mmArrow=2.5mm, line width=0.3mm, OrShift] (0,-0.1) -- ++(0,1.3);
	% ticks
%	\pic at (2.38,\Yshift) {vTick=$\frac{7+2\sqrt{6}}{5}$};
	\pic at (2,\Yshift) {vTick=$2$};
	\pic at (2.5,\Yshift) {vTick=$2.5$};
	\pic at (3,\Yshift) {vTick=$3$};
	\pic at (\Xshift,1.5) {hTick=$1.5$};
	\pic at (\Xshift,2) {hTick=$2$};
\end{tikzpicture}
%\end{document}
\caption{The graph suggests that the function $\frac{2(p+2)}{3(4-p)} - (p-2)$ does not reach the value $\frac 32$ until $p = \frac 52$, while for $\kappa \ge \frac 52$, $\kappa(p,3)$ is simply $2$. The point of intersection is $p \approx 2.8685$.}
\end{figure}
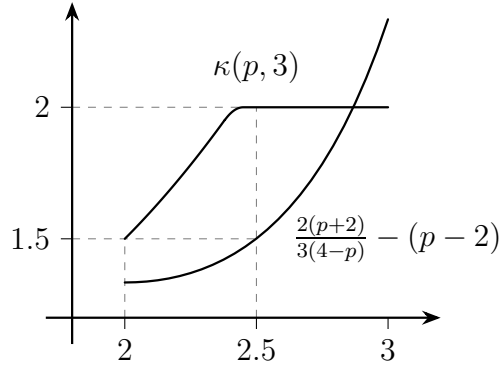
For the sake of completeness, we provide a rigorous justification. 

\textsc{Case 1:} $2 \le p \le \frac 52$. In this interval, we have 
\[
\frac{2(p+2)}{3(4-p)} - (p-2) \le \frac 32 \le \kappa(p,3).
\]
Indeed, the first inequality simplifies to $6p^2 - 23p + 20 \le 0$, which is true for $p = 2$ and $p = \frac 52$, in consequence also on the whole interval. As for the second inequality, $\kappa(p,3)$ is a nondecreasing function of $p$, and so $\kappa(p,3) \ge \kappa(2,3) = \frac 32$.

\textsc{Case 2:} $p \ge \frac 52$. Since $\frac 52 > \sqrt{6}$, we simply have $\kappa(p,3) = 2$. After rearranging, our condition is 
\[
(p-p_-)(p-p_+) \le 0,
\]
where $p_- = \frac{5-\sqrt{13}}{3} \approx 0.46$ and $p_+ = \frac{5+\sqrt{13}}{3} \approx 2.8685$. Since $\frac 52$ lies between these two roots, the condition is satisfied for $p \in [\frac 52, p_+]$, as required. 
\end{proof}

\section{Regularity for minimizing \texorpdfstring{$p$}{p}-harmonic maps \texorpdfstring{$u\colon B^3\to \S^2$}{uB3toS2}}
\label{ch:regularity-B3-S3}

The goal of this Section is to show how to improve the regularity of minimizing $p$-harmonic mappings $u\colon B^3\to S^3$ obtained in \cite{MM}. 
\begin{theorem}\label{th:regularityn=3}
Minimizing $p$-harmonic maps $u\in W^{1,p}(B^3,\S^3)$ are regular for $p\in[2,p_0]$, where $p_0\approx 2.6427$ (see \eqref{eq:maximalp} for the precise formula for $p_0$).
\end{theorem}

\begin{remark}
 Our interest in the case $p=d=3$ comes from a Cosserat type elasticity model for which the existence of nonconstant minimizing p-harmonic tangent maps is the obstruction
to full regularity; see \cite{Gastel19} for details. In \cite[Proposition 6.4]{Gastel19}, it was proven that no such tangent maps exist if $p \in [2,2.133]$ or $p\ge3$. In \cite[Theorems 4.1 and 5.1]{MM}, this was improved to $p \in [2,2.366]$ or $p\ge 2.961$. Here, we improve the first range further to $p\in [2,2.6427]$. This still leaves a gap for $p$ between $2.6427$ and $2.961$ where we cannot exclude point singularities in energy minimizing Cosserat elastic bodies.
\end{remark}

In order to prove this Theorem, we will derive improved Stability and Bochner inequalities and an appropriate mixed Kato--Cauchy--Schwarz inequality. We will be working with \emph{minimizing tangent} maps, which were defined in \Cref{de:tangentmap}.
\begin{lemma}[Improved Stability and Bochner inequalities]
\label{lem:est-in-the-limit}
Let $u \colon \S^2\to \S^3$ be a minimizing $p$-harmonic tangent map for some $2 \le p < 3$, and let $-\frac{p}{2(p-1)} < \gamma  \le 0$.
Then the following inequalities hold:
\begin{enumerate}[label=(\alph*)]
 \item\label{it:stabilitywithgamma} Stability inequality:
 \begin{equation}\label{eq:stabilitywithgamma}
\begin{split}
 3(1+\gamma)^2&\int_{\S^2} |\nabla u|^{p-2+2\gamma}|\nabla|\nabla u||^2 + (p-2)(1+\gamma)^2\int_{\S^2} |\nabla u|^{p-2+2\gamma}\abs{\left\langle \frac{\nabla u}{|\nabla u|}, \nabla |\nabla u|\right\rangle}^2\\
 & \ge (3-p)\int_{\S^2} |\nabla u|^{p+2+2\gamma} - \frac 34 (3-p)^2\int_{\S^2}|\nabla u|^{p+2\gamma};
\end{split}
\end{equation}
 \item\label{it:Bochnerwithgamma} Bochner inequality:
 \begin{equation}\label{eq:Bochnerwithgamma}
 \begin{split}
  ((p-2)&+2\gamma(p-1)) \int_{\S^2} |\nabla u|^{p-2+2\gamma}|\nabla |\nabla u||^2 + \int_{\S^2}|\nabla u|^{p-2+2\gamma}|\nabla^2 u|^2 \\
  &\le \frac12 \int_{\S^2} |\nabla u|^{p+2+2\gamma} - \int_{\S^2}|\nabla u|^{p+2\gamma}.
 \end{split}
\end{equation}
\end{enumerate}
Inequalities \ref{it:stabilitywithgamma} and \ref{it:Bochnerwithgamma} are slightly ill-defined. We adopt the (somewhat counterintuitive) convention that all ill-defined terms with a negative exponent of $|\nabla u|$ (e.g., $|\nabla u|^{p-2+2\gamma}|\nabla|\nabla u||^2$) are assumed to be zero if $|\nabla u| = 0$.
\end{lemma}

As already observed, the improved Stability and Bochner inequality involve some ill-defined expressions that are hard to work with directly. Instead, we will work with a regularized version, which roughly means that we replace $|\nabla u|$ with $|\nabla u|_\delta \coloneqq (\delta^2 + |\nabla u|^2)^{1/2}$ (for some given $\delta > 0$) and then take $\delta \to 0$. Thus, we will need the following two lemmata, whose proofs are postponed to \Cref{a:StabandBoch}.

\begin{lemma}[Regularized Stability and Bochner inequalities]\label{lem:est-with-delta}
% \Kasia{So I think $\gamma > -\frac{p}{2(p-1)}$ doesn't matters here, but later in the proof of the theorem}
Let $u$, $p$ and $\gamma$ be as in the previous Lemma. Moreover, for $\delta > 0$ denote 
$|\nabla u|_\delta \coloneqq (\delta^2 + |\nabla u|^2)^{1/2}$. Then we have
\begin{enumerate}[label=(\alph*)]
 \item\label{it:Stabilitydelta} $\delta$-Stability inequality:
\begin{equation}\label{eq:deltastability}
 \begin{split}
  &3\int_{\S^2} |\nabla u|^{p-2}|\nabla |\nabla u||^2\brac{|\nabla u|^\gamma_\delta + \gamma|\nabla u|^2|\nabla u|^{\gamma -2}_\delta}^2\\
  &\quad + (p-2)\int_{\S^2} |\nabla u|^{p-4} |\langle \nabla u, \nabla |\nabla u|\rangle|^2 \brac{|\nabla u|^\gamma_\delta
 + \gamma|\nabla u|^2|\nabla u|^{\gamma -2}_\delta}^2\\
  &\qquad \ge (3-p)\int_{\S^2} |\nabla u|^{p+2}|\nabla u|^{2\gamma}_\delta - \frac34(3-p)^2 \int_{\S^2} |\nabla u|^p|\nabla u|^{2\gamma}_\delta;
 \end{split}
\end{equation}
\item\label{it:Bochnerdelta}
$\delta$-Bochner inequality:
\begin{equation}\label{eq:Bochnerdelta}
\begin{split}
 &\int_{\S^2} |\nabla |\nabla u||^2 \brac{(1+{2\gamma})(p-2)|\nabla u|^{p-2}|\nabla u|^{2\gamma}_\delta + {2\gamma} |\nabla u|^p|\nabla u|^{{2\gamma} -2}_\delta} +|\nabla^2 u|^2 |\nabla u|^{p-2}|\nabla u|^{2\gamma}_\delta \\
 &\quad \le \int_{\S^2}- |\nabla u|^p|\nabla u|^{2\gamma}_\delta
 + \frac12|\nabla u|^{p+2}|\nabla u|^{2\gamma}_\delta.
 \end{split}
\end{equation}
\end{enumerate}
\end{lemma}

\begin{lemma}[Mixed Kato--Cauchy--Schwarz inequality]\label{le:mixedKCSwithgamma}
Let $u$, $p$ and $\gamma$ be as before. Then at all points where $\nabla u \neq 0$ we have:
\begin{equation}\label{eq:mixedKCSwithgamma}
 3|\nabla |\nabla u||^2 + (p-2) \abs{\left\langle \frac{\nabla u}{|\nabla u|}, \nabla |\nabla u|\right\rangle}^2 \\
\le C \left( |\nabla^2 u|^2 + ((p-2)+2\gamma(p-1)) |\nabla |\nabla u||^2 \right)
\end{equation}
 with $C = \frac{3}{p + 2\gamma(p-1)}$.
\end{lemma}

\begin{proof}[Proof of \Cref{lem:est-in-the-limit}]
We will show how this Lemma follows from \Cref{lem:est-with-delta}.

\underline{\Cref{lem:est-with-delta} \ref{it:Bochnerdelta} $\Rightarrow$ \Cref{lem:est-in-the-limit} \ref{it:Bochnerwithgamma}, ($\delta$-Bochner implies Bochner):}

Recall that $u$ is $C^1$ on the sphere $\S^2$. Thus, the right-hand side of \eqref{eq:Bochnerwithgamma} converges to the right-hand side of \eqref{eq:Bochnerdelta}. As for the left-hand side of \eqref{eq:Bochnerwithgamma}, let us first consider pointwise convergence of the integrand. At points where $|\nabla u| \neq 0$, this is evident. Whenever $|\nabla u| = 0$, the integrand is actually zero for each $\delta > 0$ (we have a power of $|\nabla u|$ with a positive exponent), so the limit is also zero. This agrees with our convention of interpreting \Cref{lem:est-in-the-limit}.

Note that the integrand is always nonnnegative. Indeed, due to the standard (sharp) Kato inequality $|\nabla^2 u|^2 \ge |\nabla |\nabla u||^2$ from \Cref{thm:Kato}, 
\begin{gather*}
|\nabla |\nabla u||^2 \brac{(1+{2\gamma})(p-2)|\nabla u|^{p-2}|\nabla u|^{2\gamma}_\delta + {2\gamma} |\nabla u|^p|\nabla u|^{{2\gamma} -2}_\delta} +|\nabla^2 u|^2 |\nabla u|^{p-2}|\nabla u|^{2\gamma}_\delta \\
\ge |\nabla |\nabla u||^2 |\nabla u|^{p-2} |\nabla u|^{2\gamma}_\delta
\left( (2+(1+{2\gamma})(p-2))  + {2\gamma} |\nabla u|^2 |\nabla u|^{-2}_\delta \right),
\end{gather*}
where $|\nabla u|^2 |\nabla u|^{-2}_\delta \in [0,1]$ and ${\gamma} < 0$, so it is enough to have
\[
0 \le 2+(1+{2\gamma})(p-2) + {2\gamma}
= p + {2\gamma}(p-1),
\]
which we assume from now on.
By Fatou's lemma, the integral of the limiting function is at most the lower limit of the integrals, which means that the inequality $\le$ is preserved. This finishes the proof.

\underline{\Cref{lem:est-with-delta} \ref{it:Stabilitydelta} $\Rightarrow$ \Cref{lem:est-in-the-limit} \ref{it:stabilitywithgamma}, ($\delta$-stability implies stability):}

This inequality has an opposite direction when compared to Bochner, so we cannot use Fatou's lemma. Instead, we will employ the Lebesgue dominated convergence theorem, relying on the parts of \Cref{lem:est-in-the-limit} we have already proved.

As before, the convergence of the right-hand side of \eqref{eq:stabilitywithgamma} is straightforward. The integrand of the left-hand side of \eqref{eq:stabilitywithgamma} also converges pointwise to its counterpart in \Cref{lem:est-in-the-limit} due to the presence of positive powers of $|\nabla u|$. It remains to establish a uniform bound for the integrand.

To this end, we first bound the term
\[
\left( |\nabla u|^\gamma_\delta + \gamma|\nabla u|^2|\nabla u|^{\gamma-2}_\delta \right)^2
= |\nabla u|^{2\gamma}_\delta \left( 1 + \gamma |\nabla u|^2 |\nabla u|^{-2}_\delta \right)^2 
\le |\nabla u|^{2\gamma}.
\]
This holds since $-1 < \gamma \le 0$ and $|\nabla u| \le |\nabla u|_\delta$, in particular; the expression $1 + \gamma |\nabla u|^2 |\nabla u|^{-2}_\delta$ lies between $1+\gamma$ and $1$. 
It will be convenient to denote by $\stab_\delta$ and $\stab$ the integrands appearing on the left-hand side of our stability inequalities (i.e., in ineequalities \eqref{eq:deltastability} and \eqref{eq:stabilitywithgamma}, respectively), and by $\boch$ the integrand in the limiting Bochner estimate (i.e., in \eqref{eq:Bochnerwithgamma}). The inequality above shows that $\stab_\delta \le (1+\gamma)^{-2} \stab$ pointwise. Then, the mixed Kato--Cauchy--Schwarz inequality from \Cref{le:mixedKCSwithgamma} implies $\stab \le C \cdot \boch$ with some constant $C$. To summarize, we have
\[
0 \le \stab_\delta \le (1+\gamma)^{-2} C \cdot \boch,
\]
and finally, we employ the previous step to infer that $\boch$ is indeed integrable (its integral being controlled by integrals of $|\nabla u|^{p+2\gamma}$ and $|\nabla u|^{p+2+2\gamma}$). This justifies the use of Lebesgue's dominated convergence theorem and establishes the limiting stability inequality.
\end{proof}

We are ready to proceed with the proof of the main result of the Section.

\begin{proof}[Proof of \Cref{th:regularityn=3}]
Assume that $u\colon \S^2 \to \S^3$ is a minimizing $p$-harmonic tangent map. Exactly as in, e.g., \cite[Proof of Theorem 5.1]{MM} it suffices to show that $u$ is a constant map.

Consider some $\gamma > -\frac{p}{2(p-1)}$ to be fixed later. Combining the stability inequality \eqref{eq:stabilitywithgamma} with the mixed Kato--Cauchy--Schwarz inequality \eqref{eq:mixedKCSwithgamma}, we obtain
\begin{equation}\label{eq:stabilityandmixed}
\begin{split}
 C \cdot (1+\gamma)^2
 &\brac{((p-2)+2\gamma(p-1)) \int_{\S^2} |\nabla u|^{p-2+2\gamma}|\nabla |\nabla u||^2 + \int_{\S^2}|\nabla u|^{p-2+2\gamma}|\nabla^2 u|^2 }\\
 & \ge (3-p)\int_{\S^2} |\nabla u|^{p+2+2\gamma} - (3-p)^2 \cdot \frac 34 \int_{\S^2}|\nabla u|^{p+2\gamma}
\end{split}
\end{equation}
with the constant $C = \frac{3}{p+2\gamma(p-1)}$. 
Next, we combine this result with the Bochner inequality \eqref{eq:Bochnerwithgamma}, which eliminates the second order terms: 
\begin{equation}
\begin{split}
 &(3-p)\brac{C \cdot (1+\gamma)^2}^{-1}\int_{\S^2} |\nabla u|^{p+2+2\gamma} 
 - \frac 34 (3-p)^2 \brac{C \cdot (1+\gamma)^2}^{-1} \int_{\S^2}|\nabla u|^{p+2\gamma}\\
 &\le \frac12 \int_{\S^2} |\nabla u|^{p+2+2\gamma} - \int_{\S^2}|\nabla u|^{p+2\gamma}.
 \end{split}
\end{equation}
Reordering the terms, we get
\begin{equation}
\begin{split}
A\int_{\S^2} |\nabla u|^{p+2+2\gamma} + B\int_{\S^2}|\nabla u|^{p+2\gamma} \le 0,
\end{split}
\end{equation}
where
\begin{equation}
 A=  \brac{\frac{(3-p)(p+2\gamma(p-1))}{3(1+\gamma)^2}-\frac12}, \quad B = \brac{1-\frac{(3-p)^2(p+2\gamma(p-1))}{4(1+\gamma)^2}}. 
\end{equation}
Our aim now is to choose $\gamma > -\frac{p}{2(p-1)}$ for which $A \ge 0$ and $B > 0$. If this holds, one can infer that $u$ must be constant, which completes the proof. It will be useful to consult the diagram below. 
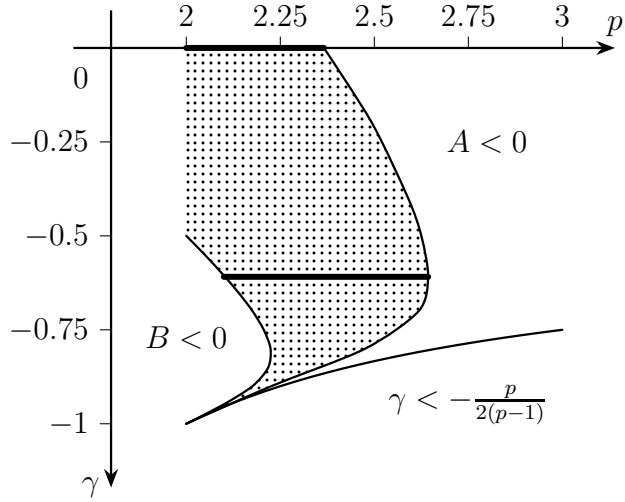
\begin{figure}[h!]
  % My favorite arrow style
\tikzstyle{mmArrow}=[-{Stealth[length=#1]}]
% style for drowing the plot line
\tikzstyle{PlotLine} = [line width=0.3mm, line cap=round]
% styles for drawing ticks on axes, with number written next to it
\tikzset{hTick/.pic={
	\draw (-0.15,0) node[left] {#1} -- (0,0);
}}
\tikzset{vTick/.pic={
	\draw (0,0.15) node[above] {#1} -- (0,0);
}}

\begin{tikzpicture}[scale=5]
	% The origin is shifted to make certain other things easier
	\def\Xshift{1.8}
	\def\Yshift{0.0}
	\tikzstyle{OrShift}=[shift={(\Xshift,\Yshift)}]
	% help lines
	\tikzstyle{Help} = [thin, gray, dashed]

  \newcommand{\Acoords}{(2.366,0) (2.494,-0.2) (2.590,-0.4) (2.625,-0.5) (2.643,-0.609) (2.618,-0.7) (2.476, -0.8) (2.226,-0.9)}
  \newcommand{\Bcoords}{(2.097,-0.95) (2.176,-0.9) (2.218,-0.85) (2.224,-0.8) (2.206,-0.75) (2.176,-0.7) (2.138,-0.65) (2.095,-0.6) (2.049,-0.55) (2,-0.5)}
  % good region
  \fill[pattern=dots]
    plot[smooth] coordinates {\Acoords (2,-1) \Bcoords}
    -- (2,0) -- cycle;
  % bad for A
  \draw[PlotLine] plot[smooth] coordinates 
    {\Acoords (2,-1)};
  \node at (2.8,-0.25) {$A < 0$};
  % bad for B
  \draw[PlotLine] plot[smooth] coordinates 
    {(2,-1) \Bcoords};
  \node at (2.0,-0.77) {$B < 0$};
  % bad for gamma
  \draw[PlotLine]
	plot[domain=2:3] (\x, {-\x/(2*\x-2)});
  \node at (2.75,-0.95) {$\gamma < -\frac{p}{2(p-1)}$};
	% axes
	\draw[mmArrow=2.5mm, line width=0.3mm, OrShift] (-0.1,0) -- ++(1.44,0) node [above] {$p$};
	\draw[mmArrow=2.5mm, line width=0.3mm, OrShift] (0,0.1) -- ++(0,-1.27) node [left] {$\gamma$};
	% ticks
	\pic at (2,\Yshift) {vTick=$2$};
	\pic at (2.25,\Yshift) {vTick=$2.25$};
	\pic at (2.5,\Yshift) {vTick=$2.5$};
	\pic at (2.75,\Yshift) {vTick=$2.75$};
	\pic at (3,\Yshift) {vTick=$3$};
	\node at ({\Xshift-0.08},-0.08) {$0$};  % tick for zero added manually, because tikz does weird things otherwise
	\pic at (\Xshift,-0.25) {hTick=$-0.25$};
	\pic at (\Xshift,-0.5) {hTick=$-0.5$};
	\pic at (\Xshift,-0.75) {hTick=$-0.75$};
	\pic at (\Xshift,-1) {hTick=$-1$};
	% finally, the ranges covered by our two choices of gamma 
	\draw[line width=0.8mm, line cap=round] 
		(2,0) -- (2.366,0)
		(2.1,-0.609) -- (2.643,-0.609);
\end{tikzpicture} 
  \caption{The graph suggests that one can find some suitable $\gamma$ for each $p \le p_0$. However, no value of $\gamma$ is suitable for all $2 \le p \le p_0$. }
\end{figure}
First, since we only consider $2 \le p < 3$, the condition $\gamma > -\frac{p}{2(p-1)}$ is satisfied for all $\gamma > -\frac 34$. If we choose $\gamma = 0$, the condition $B > 0$ is automatically satisfied, while $A \ge 0$ holds for all $p \in [2,2.36]$; for more details see our previous work \cite[Thm.~5.1]{MM}, where we followed exactly this approach. Thus, we can focus on larger values of $p$. 

Let us analyze the condition $A \ge 0$, or equivalently  
\begin{equation}
  2(3-p)(p+2\gamma(p-1)) \ge 3(1+\gamma)^2 \label{eq:pgammamax}.
\end{equation}
For a fixed $p$, this is a quadratic inequality in $\gamma$, which has a solution if and only if the discriminant is nonnegative. Thus, the maximal $p$ for which this happens can be described as a root of $2p^4 - 16p^3 + 47p^2 - 63p + 36 = 0$ or, after factoring out the trivial root $3$, of 
\[
  2p^3 - 10p^2 + 17p - 12 = 0.
\]
This cubic equation has only real root:
\begin{equation}\label{eq:maximalp}
p_0 = \frac 13 \left( 5 + \frac{\sqrt[3]{2(59+9\sqrt{43})}}{2} - \frac{1}{\sqrt[3]{2(59+9\sqrt{43})}} \right) \approx 2.6427.
\end{equation}
For the optimal $p = p_0$, the required value of $\gamma$ is
\[
\gamma = \frac{-2p_0^2+8p_0-9}{3}
= \frac 13 \left( -1 + \sqrt[3]{\sqrt{43} - 4} - \frac{3}{\sqrt[3]{\sqrt{43} - 4}} \right)
\approx -0.60874.
\]
It follows from our discussion that if we choose such $\gamma$, \eqref{eq:pgammamax} holds for all $p \in [2,p_0]$. To complete the reasoning, let us also consider the condition $B > 0$. Since $\gamma \le -\frac 12$, we have $p+2\gamma(p-1) \le 1$ and thus $B > 0$ reduces to $3-p < 2(1+\gamma)$. This clearly holds for all $p$ larger than, say, $2.3$; smaller values of $p$ have already been covered. 

In conclusion, we have shown for all $p \in [2,p_0]$ (with two different choices of $\gamma$) that any tangent map $u$ must be constant, which implies that all minimizing $p$-harmonic maps $B^3 \to \S^3$ are regular.
\end{proof}

\appendix
\section{Stability and Bochner inequalities}\label{a:StabandBoch}
\label{ch:appendix}

\begin{lemma}[Stability inequality]\label{lem:astability}
 Assume that $u \colon \S^{n-1} \to \S^d$ is a stable-stationary $p$-harmonic tangent map. 
 Then
 \begin{gather*}
\int_{\S^{n-1}} |\nabla u|^{p-4} \left( d |\nabla u|^2 |\nabla |\nabla u||^2 + (p-2) |\langle \nabla u, \nabla |\nabla u| \rangle|^2 \right) \\
\ge (d-p) \int_{\S^{n-1}} |\nabla u|^{p+2} - \frac{d(n-p)^2}{4} \int_{\S^{n-1}} |\nabla u|^p.
\end{gather*}
\end{lemma}
The proof can be found in \cite[Lem.~2,~3]{Nakauchi96}, but the result there is simplified, and thus slightly weaker. For this reason, certain technical details of the calculations will be omitted here.

\begin{proof}[Proof of \Cref{lem:astability}]
We start with the derivation from \cite[Lem.~2]{Nakauchi96}, with one final step omitted. It will be useful to consider the map  $u \colon \R^n \to \S^d$ (not restricted to the sphere), as we will exploit the stable stationary condition, i.e., the fact that the second variation is nonnegative. Following \cite[Proof of Thm. 2.7]{SU3} let us choose the perturbation family as in \cite{Xin}
\[
u_t(x) \coloneqq \frac{u(x) + t \vp(x) \alpha}{|u(x) + t  \vp(x) \alpha|},
\]
where $\alpha \in \R^{d+1}$ is a fixed unit vector and $\vp \in C_c^\infty(\R^n \setminus \{ 0 \})$ (later we will take an average over all $\alpha$ and build $\vp$ based on $|\nabla u|$).

The stability condition tells us that
\begin{align*}
0 \le \frac 1p \frac{d^2}{dt^2} \bigg|_{t=0} E_p(u_t)
= & \int_{\R^n} |\nabla u|^{p-2} |\nabla \vp|^2 (1 - |\la u, \alpha \ra|^2) \\
& + \int_{\R^n} |\nabla u|^{p-2} \vp^2 (p |\nabla \la u, \alpha \ra|^2 + (|\la u, \alpha \ra|^2-1) |\nabla u|^2) \\
& + (p-2) \int_{\R^n} |\nabla u|^{p-4} |\la \nabla \la u, \alpha \ra, \nabla \vp \ra|^2.
\end{align*}
Summing the above over all $\alpha$'s in a chosen orthonormal basis of $\R^{d+1}$, we obtain
\[
0 \le
d \int_{\R^n} |\nabla u|^{p-2} |\nabla \vp|^2
- (d-p) \int_{\R^n} |\nabla u|^{p} \vp^2
+ (p-2) \int_{\R^n} |\nabla u|^{p-4} |\la \nabla u, \nabla \vp \ra|^2.
\]
Crucially, \emph{we choose not to bound the last term via Cauchy--Schwarz inequality}. Instead, we proceed as in \cite[Lem.~3]{Nakauchi96} and choose $\vp$ in the form
\[
\vp(rz) = f(r) g(z),
\quad \text{with } f \in C_c^\infty(\R_+) \text{ and } g \in C^\infty(\S^{n-1}).
\]
Rewriting the previous inequality in this special case is straightforward; one thing worth noticing is that since $u$ is $0$-homogeneous 
% \Kasia{To jest jakas upiornie denerwujace, bo moglibysmy uzywac $u(r\xi) = \bar u(\xi)$ i wszedzie ponizej tego $\bar u$, zeby nie pisac $\nabla_\xi$... }
\begin{gather*}
|\nabla \vp(rz)|^2 = f'(r)^2 g(z)^2 + r^{-2} f(r)^2 |\nabla g(z)|^2 \\
\text{but} \quad
|\la \nabla u(rz), \nabla \vp(rz) \ra|^2
= r^{-4} f(r)^2 |\la \nabla u(z), \nabla g(z) \ra|^2,
\end{gather*}
the latter expression not involving $f'(r)^2$. Thus, we arrive at
\begin{equation*}
\begin{split}
0 &\le
\left( d \int_{\S^{n-1}} |\nabla u(z)|^{p-2} g(z)^2 \dd z \right)
\left( \int_0^\infty f'(r)^2 r^{n-p+1} \dd r \right) \\
&\quad+ \left( d \int_{\S^{n-1}} |\nabla u(z)|^{p-2} |\nabla g(z)|^2 \dd z
+ (p-2) \int_{\S^{n-1}} |\nabla u(z)|^{p-4} |\la \nabla u(z), \nabla g(z) \ra|^2 \dd z \right. \\
&\qquad\quad \left. - (d-p) \int_{\S^{n-1}} |\nabla u(z)|^p g(z)^2 \dd z \right)
\left( \int_0^\infty f(r)^2 r^{n-p-1} \dd r \right).
\end{split}
\end{equation*}
We can now optimize over all $f$. Using
\[
\inf_{f \in C_c^\infty(\R_+)} \frac{\int_0^\infty f'(r)^2 r^{n-p+1} \dd r }{\int_0^\infty f(r)^2 r^{n-p-1} \dd r}
= \frac{(n-p)^2}{4}
\]
(which can be derived by ODE considerations as in \cite[Lemma 1.3]{SU3}), we can conclude that
\begin{equation}\label{eq:stabilitywithg}
\begin{split}
&d \int_{\S^{n-1}} |\nabla u(z)|^{p-2} |\nabla g(z)|^2 \dd z
+ (p-2) \int_{\S^{n-1}} |\nabla u(z)|^{p-4} |\la \nabla u(z), \nabla g(z) \ra|^2 \dd z \\
&\quad - (d-p) \int_{\S^{n-1}} |\nabla u(z)|^p g(z)^2 \dd z
\ge -\frac{(n-p)^2}{4}
d \int_{\S^{n-1}} |\nabla u(z)|^{p-2} g(z)^2 \dd z.
\end{split}
\end{equation}
Finally, by approximating $|\nabla u|$ with smooth functions, we obtain the same inequality for $g(z) = |\nabla u(z)|$, which can be rewritten as
\begin{gather*}
d \int_{\S^{n-1}} |\nabla u(z)|^{p-2} |\nabla |\nabla u|(z)|^2 \dd z
+ (p-2) \int_{\S^{n-1}} |\nabla u(z)|^{p-4} |\la \nabla u(z), \nabla |\nabla u|(z) \ra|^2 \dd z \\
- (d-p) \int_{\S^{n-1}} |\nabla u(z)|^{p+2} \dd z
\ge -\frac{(n-p)^2}{4}
d \int_{\S^{n-1}} |\nabla u(z)|^{p} \dd z.
\end{gather*}
\end{proof}

\begin{proof}[Proof of \Cref{lem:est-with-delta} \ref{it:Stabilitydelta} ($\delta$--Stability inequality)]
 We test \eqref{eq:stabilitywithg} with the function $g=|\nabla u||\nabla u|^\gamma_\delta$, where $|\nabla u|_\delta = \brac{\delta^2 + |\nabla u|^2}^\frac12$. Then $\nabla g = \brac{|\nabla u|^\gamma_\delta + \gamma|\nabla u|^2|\nabla u|^{\gamma -2}_\delta} \nabla |\nabla u|$, which 
 leads to the desired inequality.
\end{proof}

We now turn to the Bochner inequality of \Cref{lem:est-with-delta} \ref{it:Bochnerdelta}. The starting point is the general Bochner--Weitzenb\"ock formula, which can be found e.g. in \cite[proof of Lemma 2]{Nakauchi01}, \cite[eq.~(3.13)]{EL78} or \cite[p.~123]{EelSam64} (the last reference contains a proof). For the sake of completeness, we include a direct coordinate-free derivation. 

\begin{lemma}
Consider a smooth map $u \colon \cM \to \cN$ between Riemannian manifolds $(\cM,g)$ and $(\cN,\hat{g})$. Denote the Riemannian curvature of  $\cM$ by $R$ and its Ricci curvature by $\Ric$, and similarly, by $\hat{R}$ the Riemannian curvature of $\cN$. Then, we have
\begin{equation}
\label{eq:Bochner-general}
    \frac 12 \Delta |\nabla u|^2 
    = |\nabla^2 u|^2 + \la \nabla \Delta u, \nabla u \ra 
% + \Ric_M(\nabla_a u, \nabla_a u) 
% - \Rm_{N}(\nabla_a u, \nabla_b u, \nabla_b u, \nabla_a u),
- \tensor{{\hat{R}}}{_{k\ell ij}} \nabla_a u^k \nabla_b u^\ell \nabla^b u^i \nabla^a u^j
+ \tensor{\mathrm{Ric}}{_{bc}} \nabla^b u_i \nabla^c u^i 
\end{equation}
In particular, for a map $u \colon \S^2 \to \S^3$ the following inequality holds 
\begin{equation}
\label{eq:Bochner-on-spheres}
    \frac 12 \Delta |\nabla u|^2 
    \ge |\nabla^2 u|^2 + \la \nabla \Delta u, \nabla u \ra + |\nabla u|^2 - \frac 12 |\nabla u|^4.
\end{equation}
\end{lemma}
A few words on notation are in order.
We are using Penrose's abstract index notation to ease the computation and make the curvature terms in \eqref{eq:Bochner-general} unambiguous. We will consider the pullback bundle $u^* T \cN$ with its pullback connection and curvature $\hat{R}'$. Moreover, we will use $\nabla$ to denote connections on both $T\cM$ and $T\cM \otimes u^* T\cN$, with the specific meaning being clear from the context. By $\Delta$ we denote the connection Laplacian.

Recall two basic identities for the pullback connection, taken from \cite[pp.~4-6]{EelLem83}: 
\begin{align}
    \nabla_X (du \cdot Y) - \nabla_Y (du \cdot X) & = du \cdot [X,Y], \\
    \hat{R}'(X,Y) Z & = \hat{R}(du \cdot X, du \cdot Y) Z,
\end{align}
for any tangent fields $X,Y,Z$ on $\cM$. In our notation, these read:
\begin{align}
\label{eq:P-hessian-symmetry}
    \nabla_a \nabla_b u^i & = \nabla_b \nabla_a u^i, \\
\label{eq:P-pullback-curvature}
    \tensor{{\hat{R}'}}{_{abi}^j}
    & = \tensor{{\hat{R}}}{_{k\ell i}^j} \nabla_a u^k \nabla_b u^\ell.
\end{align}

\begin{proof}
We start with computing $\frac 12 \Delta |\nabla u|^2$:  
\begin{align}
    \frac 12 \Delta |\nabla u|^2 
    & = \frac 12 \nabla_a \nabla^a (\nabla^b u_i \, \nabla_b u^i) \\
    & = \nabla_a (\nabla^b u_i \, \nabla^a \nabla_b u^i) \\
    & = \nabla_a \nabla^b u_i \, \nabla^a \nabla_b u^i + \nabla^b u_i \, \nabla_a \nabla^a \nabla_b u^i.
\end{align}
The first term we obtained is $|\nabla^2 u|^2$. In the second term, we can replace $\nabla^a \nabla_b u^i$ by $\nabla_b \nabla^a u^i$, since the Hessian of $u$ is symmetric as in \eqref{eq:P-hessian-symmetry}. If we exchange the order of differentiation once more, we obtain the desired term
\begin{equation}
    \nabla^b u_i \nabla_b \nabla_a \nabla^a u^i 
    = \la \nabla u, \nabla \Delta u \ra.
\end{equation}
However, this change results in adding a curvature term, according to \eqref{eq:P-pullback-curvature}: 
\begin{gather}
    \nabla^b u_i (\nabla_a \nabla_b \nabla^a u^i - \nabla_b \nabla_a \nabla^a u^i)
    = \nabla^b u_i \, \tensor{R}{_{abc}^a} \, \nabla^c u^i 
    + \nabla^b u_i \, \tensor{{\hat{R}'}}{_{abj}^i} \nabla^a u^j \\
    = \nabla^b u_i \, \tensor{R}{_{abc}^a} \, \nabla^c u^i 
    + \nabla^b u_i \, \tensor{{\hat{R}}}{_{k\ell j}^i} \nabla_a u^k \nabla_b u^\ell \nabla^a u^j \\
    = \tensor{\mathrm{Ric}}{_{bc}} \nabla^b u_i \nabla^c u^i 
    - \tensor{{\hat{R}}}{_{k\ell ij}} \nabla_a u^k \nabla_b u^\ell \nabla^b u^i \nabla^a u^j.
\end{gather}
This completes the derivation. 

\smallskip

In the special case of maps $u \colon \S^2 \to \S^3$, recall that $R_{abcd} = g_{ad} g_{bc} - g_{ac} g_{bd}$ and so $\tensor{\mathrm{Ric}}{_{bc}} = g_{bc}$, yielding 
% (it would be $(d-1) g_{bc}$ on a $d$-dimensional sphere). 
\begin{equation}
    \tensor{\mathrm{Ric}}{_{bc}} \nabla^b u_i \nabla^c u^i
    = |\nabla u|^2.
\end{equation}
Similarly, $\hat{R}_{k\ell ij} = \hat{g}_{kj} \hat{g}_{\ell i} - \hat{g}_{ki} \hat{g}_{\ell j}$ gives us 
\begin{equation}
%    \tensor{{\hat{R}}}{_{k\ell ij}} \nabla_a u^k \nabla_b u^\ell \nabla^a u^j \nabla^b u^i
 \tensor{{\hat{R}}}{_{k\ell ij}} \nabla_a u^k \nabla_b u^\ell \nabla^b u^i \nabla^a u^j
 = |\nabla u|^4 - \sum_{\alpha,\beta=1}^2 \la \nabla_\alpha u, \nabla_\beta u \ra^2,
\end{equation}
where the summation is over an orthonormal basis. Discarding the terms with $\alpha \neq \beta$ and applying the inequality between means of order $2$ and $4$, we get 
\begin{equation}
\sum_{\alpha,\beta=1}^2 \la \nabla_\alpha u, \nabla_\beta u \ra^2 
\ge \sum_{\alpha=1}^2 |\nabla_\alpha u|^4 \ge \frac 12 \left( \sum_{\alpha=1}^2 |\nabla_\alpha u|^2 \right)^2 = \frac 12 |\nabla u|^4.
\end{equation}
This leads to the final form of our inequality. 
\end{proof}

\begin{proof}[Proof of \Cref{lem:est-with-delta} \ref{it:Bochnerdelta} ($\delta$--Bochner inequality)]
In order to prove \eqref{eq:Bochnerdelta}, we proceed similarly as \cite[pp.~93-94]{SU3}, \cite[Lemma 1]{Nakauchi96}, or \cite[Lemma 2]{Nakauchi01}, see also \cite[p. 24]{EL78}. 
That is, we employ the Bochner--Weitzenb\"ock formula, or rather its corollary, the inequality \eqref{eq:Bochner-on-spheres} for maps $u \colon \S^2\to \S^3$. We test this with 
\[
 g=|\nabla u|^{p-2}|\nabla u|^{2\gamma}_\delta, \quad \text{ where }  |\nabla u|_\delta = (\delta^2+|\nabla u|^2)^\frac 12
\]
and integrate by parts. We do not need to exchange the order of derivatives, so the computation follows like in the flat case:  
we obtain
\begin{equation}\label{eq:Bochneraftertesting}
\begin{split}
 &-\frac 12 \int_{\S^2} \langle \nabla |\nabla u|^2,\nabla(|\nabla u|^{p-2}|\nabla u|^{{2\gamma}}_\delta)\rangle\\
 & \ge \int_{\S^2} |\nabla^2 u|^2 |\nabla u|^{p-2}|\nabla u|^{{2\gamma}}_\delta + |\nabla u|^p|\nabla u|^{2\gamma}_\delta - \frac12 |\nabla u|^{p+2}|\nabla u|^{2\gamma}_\delta -\langle \Delta u,\div (|\nabla u|^{p-2}|\nabla u|^{2\gamma}_\delta \nabla u) \rangle.
\end{split}
\end{equation}
First, we note that the term on the left-hand side of \eqref{eq:Bochneraftertesting} is
\begin{equation}\label{eq:LHSBochner}
\begin{split}
 -\frac 12 &\int_{\S^2} \langle \nabla |\nabla u|^2,\nabla(|\nabla u|^{p-2}|\nabla u|^{{2\gamma}}_\delta)\rangle\\ 
 &= -\frac 12 \int_{\S^2}  \la 2 |\nabla u||\nabla |\nabla u||\,,\,(p-2)|\nabla u|^{p-3} \nabla |\nabla u||\nabla u|^{2\gamma}_{\delta} + |\nabla u|^{p-2}{2\gamma} |\nabla u|^{{2\gamma} -2}_\delta |\nabla u| \nabla |\nabla u| \ra \\
%  &= - \int_{\S^2} |\nabla|\nabla u||^2|\nabla u|\brac{(p-2)|\nabla u|^{p-3}|\nabla u|^{2\gamma}_\delta + {2\gamma} |\nabla u|^{p-1}|\nabla u|^{{2\gamma}-2}_{\delta}}\\
 &= -\int_{\S^2} |\nabla |\nabla u||^2 \brac{(p-2) |\nabla u|^{p-2}|\nabla u|^{2\gamma}_\delta + {2\gamma} |\nabla u|^{p}|\nabla u|^{{2\gamma}-2}_\delta}.
 \end{split}
\end{equation}
Then, let us transform the last term of \eqref{eq:Bochneraftertesting} by exploiting the equation of $p$-harmonic mappings
\[
0 = \Delta_p u := \div (|\nabla u|^{p-2} \nabla u).
\]
It will also be useful to write this equation in a different form, at least at points where $\nabla u \neq 0$: 
\[
  0 = |\nabla u|^{2-p} \Delta_p u = \Delta u + (p-2) \Delta_\infty u, 
\quad \text{where }
\Delta_\infty u \coloneqq
%\left\langle \nabla^2 u, \frac{\nabla u}{|\nabla u|}\right\rangle \frac{\nabla u}{|\nabla u|} =
\big\langle \nabla |\nabla u|, \tfrac{\nabla u}{|\nabla u|}\big\rangle
\]
is the normalized (i.e., homogeneous) $\infty$-Laplace operator. This way, we have 
\[
\begin{split}
  X
  &\coloneqq -\int_{\S^2}\langle \Delta u,\div (|\nabla u|^{p-2}|\nabla u|^{2\gamma}_\delta \nabla u) \rangle\\
 &=- \int_{\S^2} \langle \Delta u, \nabla |\nabla u|^{2\gamma}_\delta \cdot |\nabla u|^{p-2}\nabla u  \rangle\\
%  &=- \int_{\S^2} \langle \Delta u, {2\gamma} |\nabla u|^{{2\gamma}-2}_\delta \langle \nabla^2 u,\nabla u\rangle \cdot |\nabla u|^{p-2}\nabla u  \rangle\\       % this particular form is correct, but the line below is more in line with our notation for \Delta_\infty
 &=- \int_{\S^2} \langle \Delta u, {2\gamma} |\nabla u|^{{2\gamma}-2}_\delta |\nabla u| \nabla |\nabla u| \cdot |\nabla u|^{p-2}\nabla u  \rangle\\
 &=- {2\gamma} \int_{\S^2} \langle \Delta u, \Delta_\infty u\rangle |\nabla u|^{{2\gamma} -2}_\delta |\nabla u|^p\\
 &= {2\gamma} (p-2) \int_{\S^2} |\Delta_\infty u|^2 |\nabla u|_\delta^{{2\gamma}-2} |\nabla u|^p.
\end{split}
\]
Since $p>2$, the pointwise inequality $|\Delta_\infty u| \le |\nabla |\nabla u||$ gives us 
\begin{equation}\label{eq:BochnerXestimate}
|X|\le 2 |\gamma|(p-2) \int_{\S^2} |\nabla u|_\delta^{{2\gamma}-2} |\nabla u|^p |\nabla |\nabla u||^2.
\end{equation}
Applying \eqref{eq:LHSBochner}, \eqref{eq:BochnerXestimate} in the inequality \eqref{eq:Bochneraftertesting}, we obtain
\[
\begin{split}
 &-\int_{\S^2} |\nabla |\nabla u||^2 \brac{(p-2)|\nabla u|^{p-2}|\nabla u|^{2\gamma}_\delta + {2\gamma} |\nabla u|^p|\nabla u|^{{2\gamma} -2}_\delta}\\
 &\ge \int_{\S^2} |\nabla^2 u|^2 |\nabla u|^{p-2}|\nabla u|^{2\gamma}_\delta + |\nabla u|^p|\nabla u|^{2\gamma}_\delta - \frac12|\nabla u|^{p+2}|\nabla u|^{2\gamma}_\delta + {2\gamma}(p-2)|\nabla|\nabla u||^2 |\nabla u|^{p-2}|\nabla u|^{2\gamma}_\delta
 \end{split}
\]
or equivalently
\[
\begin{split}
 &\int_{\S^2} |\nabla |\nabla u||^2 \brac{(1+{2\gamma})(p-2)|\nabla u|^{p-2}|\nabla u|^{2\gamma}_\delta + {2\gamma} |\nabla u|^p|\nabla u|^{{2\gamma} -2}_\delta} +|\nabla^2 u|^2 |\nabla u|^{p-2}|\nabla u|^{2\gamma}_\delta \\
 &\le \int_{\S^2}- |\nabla u|^p|\nabla u|^{2\gamma}_\delta
 + \frac12|\nabla u|^{p+2}|\nabla u|^{2\gamma}_\delta.
 \end{split}
\]
\end{proof}

\begin{proof}[Proof of \Cref{le:mixedKCSwithgamma} (mixed Kato--Cauchy--Schwarz inequality)]
We will derive an optimal constant in the inequality (at points where $\nabla u \neq 0$)
\[
3|\nabla |\nabla u||^2 + A \abs{\left\langle \frac{\nabla u}{|\nabla u|}, \nabla |\nabla u|\right\rangle}^2
\le C \left( |\nabla^2 u|^2 + B |\nabla |\nabla u||^2 \right),
\]
where $A = p-2$, $B = (p-2)+2\gamma(p-1)$, $\gamma < 0$.

We follow the proof of \cite[Lemma 5.2 (c)]{MM}. Similarly to equation (5.5) on page 3934 therein, it suffices to show that
\[
\begin{split}
3&(\theta_1 x + \theta_2 y)^2 + A \theta_1^2 (\theta_1 x + \theta_2 y)^2 \\
&\le C \Big( x^2 + 2y^2 + (x + (p-2) \theta_1 (\theta_1 x + \theta_2 y))^2 + B (\theta_1 x + \theta_2 y)^2 \Big),
\end{split}
\]
for all $\theta_1,\theta_2$ satisfying $\theta_1^2+\theta_2^2 = 1$. Taking the same steps, we introduce the rotated coordinates $z = \theta_1 x + \theta_2 y$, $w = \theta_2 x - \theta_1 y$, notice $x^2+y^2 = z^2+w^2$, $x = \theta_1 z + \theta_2 w$, and thus rearrange the above to the form
\[
3z^2 + A \theta_1^2 z^2
\le C \left( 2(z^2+w^2) + Bz^2 + (p-2)^2 \theta_1^2 z^2 + 2(p-2) \theta_1 z x \right)
\]
or, substituting $C = \frac{3}{2+B}$:
\[
A(2+B) \theta_1^2 z^2
\le 3 \left( p(p-2)\theta_1^2 z^2 + 2w^2 + 2(p-2) \theta_1 \theta_2 zw \right).
\]
Now take into account that $A = p-2$ and $2+B = p-E$, where $E = -2\gamma(p-1) \ge 0$. Then the inequality becomes
\[
0
\le (2p+E) \theta_1^2 z^2 + \frac{6}{p-2} w^2 + 6 \theta_1 \theta_2 zw,
\]
and the discriminant is
\begin{align*}
\Delta & = 36\, \theta_1^2 \theta_2^2 - 4 \cdot (2p+E) \theta_1^2  \cdot \frac{6}{p-2} \\
& = \frac{12\, \theta_1^2}{p-2} \left( 3(p-2)\theta_2^2 - 2(2p+E) \right) \\
& \le \frac{12\, \theta_1^2}{p-2} \left( 3p - 4p \right) \\
& \le 0.
\end{align*}
This shows that indeed we have
\[
3|\nabla |\nabla u||^2 + A \abs{\left\langle \frac{\nabla u}{|\nabla u|}, \nabla |\nabla u|\right\rangle}^2
\le C \left( |\nabla^2 u|^2 + B |\nabla |\nabla u||^2 \right),
\]
with $C = \frac{3}{2+B}$ in our case.
\end{proof}

\bibliographystyle{abbrv}%
\bibliography{bib}%

\end{document}